\documentclass[12pt]{article}
\usepackage{amsmath, latexsym, amsfonts, amssymb, amsthm, amscd}
\usepackage[left=3cm, right=2.5cm, top=2.5cm, bottom=2.5cm]{geometry}
\usepackage{color}

\newtheorem{proposition}{Proposition}
\newtheorem{lemma}{Lemma}

\newtheorem{theorem}{Theorem}

\newcommand{\p}{\Bbb{P}}
\newcommand{\e}{\Bbb{E}}

\newcommand{\ud}{\mathrm{d}}
\newcommand{\px}{\Bbb{P}_x}

\begin{document}

\title{Asymptotic behaviour near extinction of continuous state branching processes}
\author{G. Berzunza\footnote{ {\sc Centro de Investigaci\'on en Matem\'aticas A.C. Calle Jalisco s/n. 36240 Guanajuato, M\'exico.} E-mail: gaboh@cimat.mx}\,\, and J.C. Pardo\footnote{ {\sc Centro de Investigaci\'on en Matem\'aticas A.C. Calle Jalisco s/n. 36240 Guanajuato, M\'exico.} E-mail: jcpardo@cimat.mx}}
\maketitle
\vspace{0.2in}

\begin{abstract} 
\noindent In this note, we study the asymptotic behaviour near extinction of  (sub-) critical continuous state branching processes. In particular, we establish an analogue of Khintchin's law of the iterated logarithm near extinction time for a continuous state branching process whose branching mechanism satisfies a given condition and its reflected process at its infimum.

\bigskip

\noindent {\sc Key words and phrases}: Continuous state branching processes, Lamperti transform, L\'evy processes, conditioning to stay positive, rate of growth.

\bigskip

\noindent MSC 2000 subject classifications: 60G17, 60G51, 60G80.
\end{abstract}

\vspace{0.5cm}

\section{Introduction and main results.}
 A continuous-state branching  process (or CB-process for short)  is a non-negative valued strong Markov process with probabilities $(\mathbb{P}_x, \, x\geq 0)$ such that for any $x,y\geq 0$, $\mathbb{P}_{x+y}$ is equal in law to the convolution of $\mathbb{P}_x$ and $\mathbb{P}_y$. CB-processes may be thought of as the continuous (in time and space) analogues of classical Bienaym\'e-Galton-Watson branching processes.  Such classes of processes have  been introduced by
Jirina \cite{ji} and studied by many authors included Bingham
\cite{bi}, Grey \cite{gre}, Grimvall \cite{gri}, Lamperti \cite{la1}, to name but a few. 
More precisely, a CB-process
$Y=(Y_t, t\geq 0)$ is a Markov process taking values in
$[0,\infty]$, where $0$ and $\infty$ are two absorbing states,
and satisfying the branching property; that is to say, its
Laplace transform  satisfies

\begin{equation}
\mathbb{E}_x\Big[ e^{-\lambda  Y_t}\Big]=\exp\{-x
u_t(\lambda)\},\qquad\textrm{for }\lambda\geq 0,
\label{CBut}
\end{equation}

\noindent for some function $u_t(\lambda)$.  According to Silverstein
\cite{si}, the function $u_t(\lambda)$ is determined by the integral
equation

\begin{equation}
\int_{u_t(\lambda)}^\lambda \frac{1}{\psi(u)}{\rm d} u=t
\label{DEut}
\end{equation}

\noindent where $\psi$ satisfies the celebrate L\'evy-Khincthine formula
\begin{equation}\label{lk}
\psi(\lambda)=a\lambda+\beta
\lambda^2+\int_{(0,\infty)}\big(e^{-\lambda x}-1+\lambda x{\mathbf 1}_{\{x<1\}}\big)\Pi(\ud x),
\end{equation}
where $a\in \mathbb{R}$, $\beta\geq 0$ and $\Pi$ is a
$\sigma$-finite measure such that
$\int_{(0,\infty)}\big(1\land x^2\big)\Pi(\ud x)$ is finite.
The function $\psi$ is known as the branching mechanism of
$Y$.  

Note that the first moment of $Y_t$ can be obtained by differentiating (\ref{CBut}) with respect to  $\lambda$, i.e. $ \mathbb{E}_x[Y_t]=x e^{-\psi'(0^+)t}.$
Hence, in respective order, a CB-process is called supercritical, critical or subcritical depending on $\psi'(0^+)<0$, $\psi'(0^+)=0$ or $\psi'(0^+)>0$.  Moreover since
\[
\mathbb{P}_x\left(\lim_{t\to\infty}Y_t=0\right)=e^{-\eta x},
\]
where $\eta$ is the largest root of the branching mechanism $\psi$,  the sign of $\psi'(0+)$ yields the criterion for a.s. extinction.  More precisely a CB-process $Y$ with branching mechanism $\psi$ has a
finite time extinction almost surely if and only if
\begin{equation}\label{cond}
\int^{\infty}\frac{\ud u}{\psi(u)}<\infty\,\quad \text{ and }\quad \psi'(0+)\ge 0.
\end{equation}
 We denote by $T_{0}$  the extinction time of the CB-process $Y$, i.e. $T_{0} := \inf \{ t \geq 0: Y_{t} =0 \}. $\\

In this paper, we are interested in a detailed description of how  continuous state branching processes become extinct. Hence, in what follows we always assume that assumption (\ref{cond}) is satisfied.  One of the starting points of this paper are the results of  Kyprianou and Pardo \cite{kyp} for the CB-process in the self-similar case, i.e. when the branching mechanism is given by $\psi(\lambda)=c_+\lambda^\alpha$, for $1< \alpha\le 2$ and $c_+>0$. Note that such branching mechanism  clearly  satisfies condition (\ref{cond}).  The authors in \cite{kyp} described the upper and lower envelopes of the self-similar CB-process near extinction via integral test.  In particular they obtained the following laws of the iterated logarithm (LIL for short) for the upper envelope of $Y$ and its version reflected at its running infimum 
\[
\limsup_{t\to 0}\frac{Y_{(T_0-t)^-}}{t^{1/(\alpha-1)}\log \log(1/t)}=c_+(\alpha-1)^{1/\alpha-1}, \qquad \mathbb{P}_x-\textrm{a.s.,}
\]
and
\[
\limsup_{t\to 0}\frac{(Y-\underline{Y})_{(T_0-t)^-}}{t^{1/(\alpha-1)}\log \log(1/t)}=c_+(\alpha-1)^{1/\alpha-1}, \qquad \mathbb{P}_x-\textrm{a.s.,}
\]
where $\underline{Y}_t$ denotes the infimum of the CB-process $(Y,\mathbb{P}_x)$ over $[0,t]$.

In order to state our main results, we first introduce some basic notation. Let us  define the mapping
\[ 
\phi(t) := \int_{t}^{\infty} \frac{{\rm d}u}{\psi(u)}, \qquad \textrm{for }\quad t > 0,
\]
and note that  $\phi:(0, \infty) \rightarrow (0, \infty)$ is a bijection and thus its inverse exist, here denoted by  $\varphi$. From (\ref{DEut}),
it is straightforward to get
\[ 
u_{t}( \lambda) = \varphi(t + \phi(\lambda)) \qquad
\lambda, t >0.
 \]
Since $\phi(\infty) = 0$, we clearly have 
$u_{t}(\infty) = \varphi(t)$, for $t>0$. Hence we deduce that for every $x,t > 0$,
\[ 
\mathbb{P}_{x}(T_{0} \leq t)=\mathbb{P}_{x}( Y_{t} = 0) = \lim_{\lambda\to\infty}\e_x\Big[e^{-\lambda Y_t}\Big] =
e^{-x \varphi(t)}.
\]
Finally, let us introduce the lower and upper exponents of $\psi$ at infinity,
\[ \gamma: = \sup \left\{ c \geq 0: \lim_{\lambda \rightarrow \infty} \frac{\psi(\lambda)}{\lambda^{c}} = \infty \right\} \qquad \mbox{and} \qquad  \eta := \inf \left\{ c \geq 0: \lim_{\lambda \rightarrow \infty} \frac{\psi(\lambda)}{\lambda^{c}}  = 0 \right\}.
 \]
Since $\psi$ satisfies (\ref{lk}), we necessarily have  $1 \leq \gamma \leq
\eta \leq 2$. Now, we introduce the following exponent,
\begin{equation}\label{expod}
\delta := \sup \Big\{ c \geq 0: \exists \hspace*{1mm} Q \in (0, \infty) \hspace*{2mm} \mbox{s.t.} \hspace*{2mm} Q \psi(u)u^{-c} \leq \psi(v)v^{-c}, 1 \leq u \leq v \Big\}.
\end{equation}
Therefore, we necesarilly have $1 \leq \delta \leq \gamma \leq \eta \leq 2$. Note that   in the case where $\psi$ is regularly varying at $\infty$ with index
$\alpha >1$, then $\delta = \gamma = \eta = \alpha$.

 Our first result consist in a law of the iterated logarithm (LIL for short) at $0$ for the upper envelope of  the time-reversal process $(Y_{(T_0 -t)^{-}}, 0\leq t\leq T_0),$ under $\px$. 
\begin{theorem}\label{teo2}
Assume that $\delta > 1$, then
\[ \limsup_{t \rightarrow 0} \frac{Y_{(T_{0}-t)^{-}} \varphi(t)}{\log \log
\varphi(t)} = 1 \hspace*{7mm} \px\mbox{-a.s.,} \]
\noindent for every $x > 0$.
\end{theorem}

Recall that  $\underline{Y}$ and $((Y-\underline{Y})_{(T_{0}-t)^{-}}, 0 \leq t < T_{0})$ denote the running infimum of $Y$ and 
the time-reversal process of $Y$ reflected at its running minimum, respectively. 
 We also introduce the so-called scale function $W: [0,\infty)\mapsto [0,\infty)$ which is the unique absolutely continuous increasing function whose Laplace transform is $1/\psi$.  
Let us suppose that for all $\beta < 1$, the scale function $W$ satisfies the following hypothesis
\[ \hspace{-4cm}(\mathbf{H})\quad\qquad\qquad\qquad \limsup_{x \rightarrow 0} \frac{W(\beta
x)}{W(x)} < 1. 
\]
We remark that the above hypothesis is satisfied, in particular, when $\psi$ is regularly varying at $\infty$.
\begin{theorem} \label{teo3}
  Suppose that $\delta > 1$, then under the hypothesis $(\mathbf{H})$, we have 
  \[ \limsup_{t \rightarrow 0} \frac{(Y -\underline{Y} )_{(T_{0}-t)^{-}}\varphi(t)}{\log \log
\varphi(t)} = 1 \hspace*{7mm} \px\mbox{-a.s.,} \]
\noindent for every $x > 0$.
\end{theorem}
It is important to note that in the self-similar case, i.e. when $\psi(\lambda) = c_{+} \lambda^{\alpha}$ for $1<\alpha\le 2$ and $c_{+}
 >0$, we have
 \[  \frac{\log \log \varphi(t)}{\varphi(t)}= \frac{\log \log
 (c_{+}(\alpha-1)t)^{-\frac{1}{\alpha-1}}}{(c_{+}(\alpha-1)t)^{-\frac{1}{\alpha-1}}}  \qquad \mbox{for} \quad t >0. \]
We can replace the  above function by 
 \[ \frac{\log \log \frac{1}{t}}{ \varphi(t)} = (c_{+}(\alpha-1)t)^{\frac{1}{\alpha-1}}
  \log \log{\frac{1}{t}}\qquad\mbox{for}\quad t >0, \]
since they are asymptotically equivalent at 0. Similarly, we can take the previous function in the regularly varying case.

\section{Proofs}\label{second}

Let $(\mathbf{P}_x,\, x\in \mathbb{R})$ be a family of probability measures on the
space of c\`adl\`ag mappings from $[0,\infty)$ to $\mathbb{R}$, denoted by
$\mathcal{D}$, such that for each $x\in \mathbb{R}$, the canonical process
$X$ is a L\'evy process with no negative jumps issued from $x$. Set
$\mathbf{P}:=\mathbf{P}_0$, so $\mathbf{P}_x$ is the law of $X+x$ under $\mathbf{P}$. The Laplace
exponent $\psi:[0,\infty) \to (-\infty,\infty)$ of $X$ is specified
by $\mathbf{E}[e^{-\lambda X_t}]=e^{t\psi(\lambda)}$, for $ \lambda\in \mathbb{R}$,
and can be expressed in terms of the L\'evy-Khintchine formula  (\ref{lk}).

Henceforth, we shall assume that $(X,\mathbf{P})$ is not a subordinator
(recall that  a subordinator is a L\'evy process with increasing
sample paths). In that case, it is known that the Laplace exponent
$\psi$ is strictly convex and tends to $\infty$ as $\lambda$ goes to
$\infty$. In this case, we define for $q\geq 0$
\[
\Phi(q)=\inf\big\{\lambda\geq 0: \psi(\lambda)>q\big\}
\]
the right-continuous inverse of $\psi$ and then $\Phi(0)$ is the
largest root of the equation $\psi(\lambda)=0$. Theorem VII.1 in \cite{Be} implies that condition $\Phi(0)>0$ holds if and only if the process drifts to $\infty$. Moreover, almost surely, the paths of $X$ drift to $\infty$, oscillate or drift to $-\infty$ accordingly as $\psi'(0+)<0$, $\psi'(0+)=0$ or $\psi'(0+)>0$.

Lamperti \cite{la1} observed that  continuous state branching
processes are connected to L\'evy processes with no negative jumps
by a simple time-change. More precisely, consider a spectrally
positive L\'evy process $(X,\mathbf{P}_x)$ started at $x>0$ and with Laplace
exponent $\psi$. Now, we introduce the clock
\[
A_{t}=\int_0^{t}\frac{\ud s}{X_s} \qquad t\in [0,\tau_0),
\]
where $\tau_0=\inf\{t\ge 0: X_t\le 0\}$,  and its right-continuous inverse
\begin{equation*}
\theta(t)=\inf\{s\geq 0:\,A_s>t\}.
\end{equation*}
Then, the time changed process $Y=(X_{\theta(t)}, t\geq 0)$, under
$\mathbf{P}_x$, is a continuous state branching process with initial
population of size $x$. The transformation described above will
henceforth be referred to as the {\it CB-Lamperti representation}.



Now, define $\widehat{X}:=-X$, the dual process of $X$. Denote by
$\widehat{\mathbf{P}}_x$ the law of $\widehat{X}$ when issued from $x$ so
that $(X,\widehat{\mathbf{P}}_x) = (\widehat{X}, \mathbf{P}_{-x})$. The dual
process conditioned  to stay positive  is a Doob $h$-transform of $(X, \widehat{\mathbf{P}}_x)$ killed when it
first exists $(0,\infty)$ with the harmonic function $W$. In this case, assuming  that  $\psi'(0+) \geq 0$, one has
\[
\widehat{\mathbf{P}}_x^{\uparrow}(X_{t}\in {\rm d}y)=
\frac{W(y)}{W(x)}\widehat{\mathbf{P}}_x(X_t\in {\rm d}y\,,t<\tau_0)
\,,\;\;\;t\ge0,\,\;\;\;x,y>0.
\]
Under $\widehat{\mathbf{P}}_x^\uparrow$, $X$ is a process taking values in $(0,\infty)$. It will be referred to as the dual L\'evy process started at $x$ and conditioned to stay positive. The measure
$\widehat{\mathbf{P}}_x^\uparrow$ is always a probability measure and there
is always weak convergence as $x\downarrow 0$ to a probability
measure which we denote by $\widehat{\mathbf{P}}^\uparrow$.

 Let
\[
 \sigma_x = \sup\{t>0 : X_t \leq x\},
 \]
be the last passage time of $X$ below  $x\in \mathbb{R}$. The
next proposition, whose proof may be found in \cite{kyp},  gives us a time-reversal property from
extinction for the C.B. process (see  Theorem
1). Recall that, under $\widehat{\mathbf{P}}^\uparrow$, the canonical process
$X$ drifts towards $\infty$ and also that $X_t>0$ for $t>0$.

\begin{proposition}\label{CBreversedentrance} If  condition (\ref{cond}) holds, then for each $x>0$
\begin{equation*}
\Big\{(Y_{(T_0 -t)^{-}}  : 0 \leq t< T_0 ), \mathbb{P}_x\Big\} \stackrel{d}{=}
\Big\{(X_{\theta(t)}, 0\leq t<A_{\sigma_x}),\widehat{\mathbf{P}}^\uparrow\Big\}.
\end{equation*}
\end{proposition}

Recall that $\psi$ is a branching mechanism satisfying conditions 
(\ref{lk}) and  (\ref{cond}); and  whose exponent $\delta$ defined
by (\ref{expod}) is strictly larger than $1$. Therefore, since
$\delta >1$, there exist $c \in (1, \infty)$ and $C \in (0, \infty)$
such that
\begin{eqnarray} \label{desimp}
  \psi(\lambda) \leq C \psi(\lambda u) u ^{-c},
\end{eqnarray}
\noindent for any $u$, $\lambda \in [1, \infty)$. Now, we introduce
the so-called first and last passage times of the  process $\{(X_{\theta(t)}, t\ge 0),\widehat{\mathbf{P}}^\uparrow\} $ by
\[ S_{y} = \inf \Big\{ t \geq 0:\, X_{\theta(t)} \geq y \Big\} \qquad
\mbox{and} \qquad U_{y} = \sup \Big\{ t \geq 0: \,X_{\theta(t)} \leq y \Big\}, \]
 for $y\ge 0$. Note that the processes $(S_y, y\ge 0)$ and $(U_y, y\ge 0)$  are
increasing with independent increments since the process $\{(X_{\theta(t)}, t\ge 0),\widehat{\mathbf{P}}^\uparrow\} $ has no positive jumps. 

Observe that  $A_{\sigma_{x}}$ and $U_x$ are equal and that the latter, under  $\widehat{\mathbf{P}}^\uparrow$, has the same law as   $T_{0}$ under $\p_{x}$. This clearly implies that the distribution of
${U}_{x}$, under  $\widehat{\mathbf{P}}^\uparrow$, satisfies
\begin{eqnarray} \label{distlast}
\widehat{\mathbf{P}}^\uparrow \Big({U}_{x} \leq t \Big) = e^{-x \varphi(t)},
\end{eqnarray}
\noindent for every $x,t > 0$. 

Let $(J_t, t\ge 0)$ be the future infimum process of $(X_{\theta(t)}, t\ge 0)$ which is defined as follows
\[
J_t=\inf_{s\ge t} X_{\theta(s)},\qquad \textrm{for } \quad t\ge 0.
\]
\begin{proposition}Assume that $\delta > 1$, then
\[ \limsup_{t \rightarrow 0} \frac{J_t \varphi(t)}{\log \log
\varphi(t)} = 1, \hspace*{7mm} \widehat{\mathbf{P}}^\uparrow\mbox{-a.s.} \]
\end{proposition}
\noindent \emph{Proof:} For simplicity, we let 
\[ f(t) = \frac{\log \log \varphi(t)}{\varphi(t)} \qquad \mbox{for} \quad t >0. \]
 In order to prove this result, we need the following
two technical lemmas. Recall that $\phi$ is the inverse function of $\varphi$.
\begin{lemma} \label{lemma2} For every integer $n \geq 1$ and $r>1$, put
  \[ t_{n} = \phi(r^{n}) \qquad\textrm{and}\qquad  a_{n} = f(t_{n}). \]
\begin{itemize}
  \item[(i)] The sequence $(t_{n}: n \geq 1)$ decreases.
  \item[(ii)] The series $\sum_{n}  \widehat{\mathbf{P}}^\uparrow\big( J_{t_n} > r
  a_{n}\big)$ converges.
\end{itemize}
\end{lemma}
\noindent \emph{Proof of Lemma \ref{lemma2}:}  The first assertion follows readily from the fact that $\phi$ is
decreasing. In order to prove (ii), we note that  for $n \geq
1$, $\varphi(t_{n}) = \varphi(\phi(r^{n})) = r^{n}.$ This entails
\begin{eqnarray} \label{ec2lema2}
\log \log \varphi(t_{n}) = \log \log r^{n}.
\end{eqnarray}
Now, since the last passage times process is the right continuos
inverse of the future infimum of $(X_{\theta(t)}, t\ge 0)$, we have 
\[ \widehat{\mathbf{P}}^\uparrow\Big(J_{t_n}> r a_{n} \Big) = \widehat{\mathbf{P}}^\uparrow\Big( {U}_{ ra_{n}} < t_{n}\Big). \]
Hence (\ref{distlast}) and  (\ref{ec2lema2}) imply
\begin{eqnarray*}
 \sum_{n} \widehat{\mathbf{P}}^\uparrow\Big( {U}_{ra_{n}} < t_{n}\Big) & \leq & \sum_{n} \left( \frac{1}{n \log r} \right)^{r},
\end{eqnarray*}
\noindent which converges, and our statement follows.
\hfill $\Box$

\begin{lemma} \label{lemma3}
 For every integer $n \geq 2$ and $r>1$, put
\[ s_{n} = \phi(e^{n^{r} })\qquad \textrm{and}\qquad b_{n} = f(s_{n}). \]
We have that the series $\sum_{n}  \widehat{\mathbf{P}}^\uparrow\Big({U}_{b_{n}/r} \leq s_{n}\Big)$ diverges.
\end{lemma}
\noindent \emph{Proof of Lemma \ref{lemma3}:} First we note that  for $n \geq 1,$ $ \varphi(s_{n}) = \varphi(\phi(e^{n^{r}})) = e^{n^{r}}.$
This entails  $\log \log \varphi(s_{n}) = \log n^{r}.$ Hence, the identity  (\ref{distlast})  implies
\begin{eqnarray*}
 \sum_{n} \widehat{\mathbf{P}}^\uparrow\Big( {U}_{b_{n}/r} \leq s_{n}\Big) & = & \sum_{n} \frac{1}{n} ,
\end{eqnarray*}
\noindent which diverges, and our claim follows.
\hfill $\Box$ 

We are now able to establish the law of the iterated logarithm. In order to
prove the upper bound, we use  Lemma \ref{lemma2}.
Take any $t \in [t_{n+1}, t_{n}]$, so, provided that $n$ is large
enough
\[
 f(t) \geq \frac{\log \log \varphi(t_{n})}{\varphi(t_{n+1})}, 
 \]
since $\varphi$ decreases. Note that the
denumerator is equal to $r^{n+1}$ and  the numerator is equal to
$\log \log r^{n}$. We thus have
\begin{equation}\label{foverf}
\limsup_{t \rightarrow 0} \frac{f(t_{n})}{f(t)} \leq r. 
\end{equation}
On the other hand, an  application of the Borel-Cantelli Lemma to Lemma
\ref{lemma2} shows that
\[ \limsup_{n \rightarrow \infty}
\frac{J_{t_n}}{f(t_{n})} \leq r ,\qquad  \widehat{\mathbf{P}}^\uparrow\mbox{-a.s.,} \]
 and we deduce that
\[ \limsup_{t \rightarrow 0} \frac{J_t}{f(t)} \leq \left(\limsup_{n \rightarrow
\infty} \frac{J_{t_n}}{f(t_{n})}\right) \left(\limsup_{t \rightarrow 0} \frac{f(t_{n})}{f(t)}\right) \leq r^{2}, \qquad  \widehat{\mathbf{P}}^\uparrow\mbox{-a.s.} \]

To prove the lower bound, we use  Lemma \ref{lemma3}
and observe that the sequence $(b_{n}, n \geq 2)$ decreases. First, from Lemma \ref{lemma3}, we have
\[ \sum_{n} \widehat{\mathbf{P}}^\uparrow\Big( U_{b_{n}/r} - U_{b_{n+1}/r} \leq s_{n}\Big) \geq \sum_{n}
\widehat{\mathbf{P}}^\uparrow\Big( U_{b_{n}/r} \leq s_{n}\Big) = \infty, \]
so by the Borel-Cantelli Lemma for independent events, we obtain
\[ \liminf_{n \rightarrow \infty} \frac{U_{b_{n}/r} -
U_{b_{n+1}/r}}{s_{n}} \leq 1,\qquad\widehat{\mathbf{P}}^\uparrow\mbox{-a.s.} \]
If we admit for a while that
\begin{eqnarray} \label{ec1teo2}
\limsup_{n \rightarrow \infty} \frac{U_{b_{n+1}/r}}{s_{n}}
= 0, \qquad\widehat{\mathbf{P}}^\uparrow\mbox{-a.s.,}
\end{eqnarray}
we can conclude 
\[ \liminf_{n \rightarrow \infty} \frac{U_{b_{n}/r}}{s_{n}} \leq 1, \qquad \widehat{\mathbf{P}}^\uparrow\mbox{-a.s.} \]
This implies that the set $\{s: U_{f(s)/r} \leq s \}$ is unbounded
$\widehat{\mathbf{P}}^\uparrow$-a.s. Plainly, the same then holds for $\{s: J_s \geq f(s)/r
\}$, and as a consequence
\begin{eqnarray} \label{ec2teo2}
 \limsup_{t \rightarrow 0} \frac{J_t}{f(t)} \geq
\frac{1}{r}, \qquad\widehat{\mathbf{P}}^\uparrow\mbox{-a.s.}
\end{eqnarray}
Next we establish the behaviour in (\ref{ec1teo2}). First, since $\delta > 1$ we
observed from  inequality (\ref{desimp}) that there exist $c \in (1,
\infty)$ and $C \in (0,\infty)$ such that
\begin{equation}
  \varepsilon \phi(t)  =  \varepsilon \int_{t}^{\infty} \frac{du}{\psi(u)} \geq  \left(\frac{C}{\varepsilon}\right)^{\frac{1}{c-1}} \int_{t}^{\infty} \frac{du}{\psi\big( \left(\frac{C}{\varepsilon}\right)^{\frac{1}{c-1}} u \big)} =  \int_{ (\frac{C}{\varepsilon})^{\frac{1}{c-1}} t}^{\infty} \frac{du}{\psi(u)} =  \phi\big((C/\varepsilon)^{\frac{1}{c-1}} t\big), \label{ec3teo2}
\end{equation}
 for any $t \geq 1$ and $0 < \varepsilon < \min(C, 1)$. On the other hand,  we
see that for $n$ large enough and  $0 < \varepsilon < \min(C, 1)$
\[
 \widehat{\mathbf{P}}^\uparrow\Big( {U}_{b_{n+1}/r} > \varepsilon s_{n} \Big) =  1 - \exp \left\{- \frac{b_{n+1}}{r} \varphi(\varepsilon s_{n}) \right\}  \leq  \frac{b_{n+1}}{2} \varphi(\varepsilon s_{n}). \\
\]
Hence from inequality (\ref{ec3teo2}), we have
\[
 \widehat{\mathbf{P}}^\uparrow\Big( {U}_{b_{n+1}/r} > \varepsilon s_{n} \Big)\leq  \left( \frac{C}{\varepsilon} \right)^{\frac{1}{c -1}} \exp \Big\{ n^{r} - (n+1)^{r} \Big\} \log (n+1)\le  \left( \frac{C}{\varepsilon}
\right)^{\frac{1}{c -1}} e^{-(r-1)n^{r}+1},
\]
where the las identity follows since $(n+1)^{r}-n^{r} \geq rn^{r-1}$. In conclusion, we have
that  the series $\sum_{n} \widehat{\mathbf{P}}^\uparrow\big( U_{b_{n+1}/r} > s_{n} \big)$
converges, and according to the Borel-Cantelli lemma,
\[ \limsup_{n \rightarrow \infty} \frac{U_{b_{n+1}/r}}{s_{n}}
\leq \varepsilon, \qquad \widehat{\mathbf{P}}^\uparrow\mbox{-a.s.,} \]

\noindent which establishes (\ref{ec1teo2}) since $0 < \varepsilon < \min(C, 1)$ can be chosen arbitrarily small. The proof of (\ref{ec2teo2}) is now complete. The two preceding bounds show that
\[ \frac{1}{r} \leq \limsup_{t \rightarrow 0} \frac{J_t}{f(t)} \leq r^{2}, \qquad\widehat{\mathbf{P}}^\uparrow\mbox{-a.s.} \]
Hence the result follows taking  $r$ close enough to $1$.
\hfill $\Box$ \\
\noindent \emph{Proof of Theorem \ref{teo2}:} In order to establish our result, we first prove the following law of the iterated logarithm holds
\begin{equation}\label{lilX}
\limsup_{t \rightarrow 0}
\frac{X_{\theta(t)}}{ f(t)}=1, \qquad \widehat{\mathbf{P}}^\uparrow-\textrm{a.s.},
\end{equation}
and then use Proposition 1.

The lower bound of (\ref{lilX}) is easy to deduce from Proposition 2. More precisely, 
\[ 1 = \limsup_{t \rightarrow 0}
\frac{J_t}{ f(t)} \leq
\limsup_{t \rightarrow 0} \frac{X_{\theta(t)}}{ f(t)}\qquad \widehat{\mathbf{P}}^\uparrow-\mbox{a.s.} \]
Now, we prove the upper bound. Let $r > 1$ and denote by  $\overline{X}_{\theta(\cdot)}$ for
the supremum process of $\{(X_{\theta(s)},s\ge0),\widehat{\mathbf{P}}^{\uparrow}\}$, which is defined by $\overline{X}_{\theta(t)} =
\sup_{0 \leq s \leq t} X_{\theta(s)}$ for $t \geq 0$. 
\begin{lemma} \label{lemma4}
  Let $c^\prime>1$, $0 < \epsilon < 1-1/c^{\prime}$ and $r > 1$. For $t_{n} = \phi(r^{n})$,  $n\ge 1$,  then there exist a
  positive real number $K$ such that
  \[
\widehat{\mathbf{P}}^\uparrow\Big( J_{t_n}  > (1 - \epsilon)c^\prime  f(t_{n})\Big) \geq   \epsilon K^{2} \widehat{\mathbf{P}}^\uparrow\Big(\overline{X}_{\theta(t_n)} >c^\prime f(t_n)\Big). 
 \]
\end{lemma}
\noindent \emph{Proof:}  From the Markov property, we have
\begin{equation*}
\begin{split} \label{ec1lema4}
\widehat{\mathbf{P}}^\uparrow\bigg(J_{t_n}   > (&1 - \epsilon)c^\prime   f(t_{n})\bigg)  \geq \widehat{\mathbf{P}}^\uparrow\Big(\overline{X}_{\theta(t_n)} >  c^\prime   f(t_{n}), J_{t_n}  > (1 - \epsilon)  c^\prime f(t_{n})\Big) \\
                                                          & =  \int_{0}^{t_{n}}\widehat{\mathbf{P}}^\uparrow\Big(S_{c^\prime  f(t_{n})} \in \mathrm{d} t\Big) \widehat{\mathbf{P}}^\uparrow_{c^\prime f(t_n)}\Big(J_{t_n}  > (1 - \epsilon)c^\prime   f(t_{n})\Big)  \\
                                                          &\geq \widehat{\mathbf{P}}^\uparrow\Big(\overline{X}_{\theta(t_n)} >c^\prime f(t_n)\Big) \widehat{\mathbf{P}}^\uparrow_{c^\prime f(t_n)}\left( \inf_{0\le s} X_{\theta(s)}  > (1 - \epsilon)c^\prime   f(t_{n})\right).
                                                          \end{split}
\end{equation*}
Now from the Lamperti transform and   Lemma VII.12 in \cite{Be}, we have
\begin{equation*}
\begin{split}
\widehat{\mathbf{P}}^\uparrow_{c^\prime f(t_n)}\left(\inf_{0\le s} X_{\theta(s)}  > (1 - \epsilon) c^\prime   f(t_{n})\right)  = \widehat{\mathbf{P}}^\uparrow_{c^\prime f(t_n)}\Big(\tau_{[0,(1-\epsilon)c^\prime f(t_n))}=\infty\Big) = \frac{W\big(c^\prime\epsilon   f(t_{n}) \big)}{W\big(c^\prime f(t_{n})\big)},
  \end{split}
\end{equation*}
where $\tau_{[0,z)}=\inf\{t\ge 0: X_t\in [0,z)\}$.
 
On the other hand an application of Proposition III.1 in \cite{Be}
gives that there exist a positive real number $K$ such that
\begin{equation}\label{eqlemma3}
K \frac{1}{x \psi(1/x)} \leq W(x) \leq K^{-1} \frac{1}{x
\psi(1/x)}, \qquad \mbox{for all}\quad x >0,
 \end{equation}
\noindent then it is clear 
\[ \frac{W\big(c^\prime \epsilon   f(t_{n}) \big)}{W\big(c^\prime   f(t_{n})\big)} \geq
K^{2} \epsilon^{-1} \frac{\psi\big(1/c^\prime  f(t_{n})\big)}{ \psi\big(\epsilon^{-1}/c^\prime   f(t_{n})\big)}. \]
From the above  inequality and Lemma 3 in \cite{Pardo}, we deduce
  \[
\widehat{\mathbf{P}}^\uparrow_{c^\prime f(t_n)}\left(\inf_{0\le s} X_{\theta(s)}  > (1 - \epsilon)c^\prime  f(t_{n})\right)   \geq   \epsilon K^{2},
 \]
which clearly implies our result.
\hfill $\Box$

Now, we prove the upper bound for the LIL of $((X_{\theta(t)}, t\ge 0),\widehat{\mathbf{P}}^\uparrow)$. Let $c^\prime>1$ and fix $0 <
\epsilon < 1-1/c^\prime$.  Recall from  (\ref{distlast}) that
\[
\begin{split}
\widehat{\mathbf{P}}^\uparrow\Big( J_{t_n}> (1-\epsilon)c^\prime  f(t_n) \Big)& = \widehat{\mathbf{P}}^\uparrow\Big( U_{(1- \epsilon)c^\prime  f(t_{n}) }<
t_{n-1}\Big)=(n \log r )^{-(1-\epsilon)c^\prime }.
\end{split}
\]
Hence from  Lemma \ref{lemma4}, we deduce
\[ 
\sum_{n\ge 1}\widehat{\mathbf{P}}^\uparrow\Big(S_{t_{n}} >c^\prime f(t_n)\Big)\leq K^{-2} \epsilon^{-1} \sum_{n} (n \log r
)^{-(1-\epsilon)c^\prime }<\infty, 
\]
since $(1-\epsilon)c^\prime>1$. Therefore an application of the Borel-Cantelli Lemma  shows 
\[ 
\limsup_{n \rightarrow \infty} \frac{\overline{X}_{\theta(t_n)}  }{f(t_{n})} \leq  c^\prime,
\hspace*{5mm} \widehat{\mathbf{P}}^\uparrow\mbox{-a.s.} 
\]
Using (\ref{foverf}), we deduce 
\[ \limsup_{t \rightarrow 0} \frac{\overline{X}_{\theta(t)} }{f(t)} \leq \left(\limsup_{n \rightarrow
\infty} \frac{\overline{X}_{\theta(t_n) }  }{f(t_{n})}\right) \left(\limsup_{t \rightarrow 0} \frac{f(t_{n})}{f(t)}\right) \leq r c^\prime, \quad\widehat{\mathbf{P}}^\uparrow\mbox{-a.s.} \]
Hence from  the first part of the proof  and  taking  $r$ close enought to $1$  above, we deduce 
\[ \limsup_{t \rightarrow 0} \frac{X_{\theta(t)} }{f(t)} \in [1,c^\prime], \quad \widehat{\mathbf{P}}^\uparrow\mbox{-a.s.}, 
\]
By the Blumenthal zero-one law, it must be a constant number $k$, $\widehat{\mathbf{P}}^\uparrow$-a.s.

Now we prove that the constant $k$ equals 1. Fix $\epsilon \in (0, 1/2)$ and define
\begin{eqnarray*}
 R_{n} = \inf \left \{ \frac{1}{n} \leq s: \frac{J_s}{k f(s)} \geq (1- \epsilon) \right \}.
\end{eqnarray*}
Note that for $n$ sufficiently large $1/n < R_{n} < \infty$ and that  $R_{n}$
converge to $0$, $\widehat{\mathbf{P}}^\uparrow$-a.s.,  as $n$ goes to $\infty$.  From Lemma VII.12 in \cite{Be}, the strong Markov property and since the process $\{ (X_{\theta(t)}, t\ge 0) , \widehat{\mathbf{P}}^\uparrow \}$ has no positive jumps, we have 
\[
\begin{split}
 \widehat{\mathbf{P}}^\uparrow \left( \frac{J_{R_{n}}}{kf(R_{n})} \geq (1 - 2\epsilon) \right) & =  \widehat{\mathbf{P}}^\uparrow\left(  J_{R_{n}} \geq \frac{(1 - 2\epsilon)X_{\theta(R_{n}) }}{(1- \epsilon) }  \right) \\
 & = \widehat{\mathbf{E}}^\uparrow \left(\widehat{\mathbf{P}}^\uparrow\left( J_{R_{n}} \geq \frac{(1 - 2\epsilon)X_{\theta(R_{n})}}{(1- \epsilon) }   \Big\vert X_{\theta(R_{n})}\right) \right) \\
                                                                                           & =  \widehat{\mathbf{E}}^\uparrow \left( \frac{W(\ell(\epsilon) X_{\theta(R_{n})})}{W( X_{\theta(R_{n})})} \right), \\
\end{split}
\]
\noindent where $\ell(\epsilon) =  \epsilon/(1- \epsilon)$. Applying (\ref{eqlemma3}) and Lemma 3 in \cite{Pardo}, give us
\[
\widehat{\mathbf{E}}^\uparrow \left( \frac{W(\ell(\epsilon) X_{\theta(R_{n})})}{W( X_{\theta(R_{n})})} \right)\ge K^2\ell(\epsilon),
\]
wich implies 
\[
\lim_{n\to\infty} \widehat{\mathbf{P}}^\uparrow \left( \frac{J_{R_{n}}}{kf(R_{n})} \geq (1 - 2\epsilon) \right)>0. 
\]
Since $R_n\ge 1/n$,
\[
 \widehat{\mathbf{P}}^\uparrow \left( \frac{J_{t}}{kf(t)} \geq (1 - 2\epsilon), \,\,\textrm{for some }t\ge 1/n \right)\ge\widehat{\mathbf{P}}^\uparrow \left( \frac{J_{R_{n}}}{kf(R_{n})} \geq (1 - 2\epsilon) \right).
\]
Therefore, for all $\epsilon\in (0,1/2)$
\[
 \widehat{\mathbf{P}}^\uparrow \left( \frac{J_{t}}{kf(t)} \geq (1 - 2\epsilon), \,\,\textrm{i.o., as  }t\to 0\right)\ge\lim_{n\to\infty}\widehat{\mathbf{P}}^\uparrow \left( \frac{J_{R_{n}}}{kf(R_{n})} \geq (1 - 2\epsilon) \right)>0.
\]
The event on the left hand side is in the lower-tail sigma-field of $(X,\widehat{\mathbf{P}}^\uparrow)$ which is trivial from Bertoin's contruction (see for instance Section 8.5.2 in \cite{Don}). Hence
\[
\limsup_{t\to 0}\frac{J_t}{f(t)}\ge k(1-2\epsilon), \qquad \widehat{\mathbf{P}}^\uparrow\textrm{-a.s.,}
\]
and since $\epsilon$ can be chosen arbitraily small, we deduce that $1\ge k$.
\hfill $\Box$

\noindent \emph{Proof of Theorem \ref{teo3}:} Here we follow similar arguments as those used in the last part of the previous result.
Assume that the  hypothesis $(\mathbf{H})$ is satisfied. From Theorem \ref{teo2}, it is clear that
\[
 \limsup_{t \rightarrow 0} \frac{X_{\theta(t)} - J_{t}}{ f(t)} \leq
\limsup_{t \rightarrow 0} \frac{X_{\theta(t)}}{ f(t)} = 1,\qquad \widehat{\mathbf{P}}^\uparrow\textrm{-a.s.} 
\]
Fix $\varepsilon \in (0, 1/2)$ and define
\begin{eqnarray*}
 R_{n} = \inf \left \{ \frac{1}{n} \leq s: \frac{X_{\theta(s)}}{f(s)} \geq (1- \epsilon) \right \}.
\end{eqnarray*}
First note that for $n$ sufficiently large $1/n < R_{n} < \infty$
$\widehat{\mathbf{P}}^\uparrow$-a.s. Moreover, from Theorem \ref{teo2} we have that $R_{n}$
converge to $0$ as $n$ goes to $\infty$, $\widehat{\mathbf{P}}^\uparrow$-a.s. 

From Lemma VII.12 in \cite{Be}, then strong Markov property and since $\{ (X_{\theta(t)}, t\ge 0) , \widehat{\mathbf{P}}^\uparrow \}$ has no positive jumps, we have 
\[
\begin{split}
\widehat{\mathbf{P}}^\uparrow \left( \frac{X_{\theta(R_n)}-J_{R_{n}}}{f(R_{n})} \geq (1 - 2\epsilon) \right) & =  \widehat{\mathbf{P}}^\uparrow \left( J_{R_{n}} \leq \frac{\epsilon}{1- \epsilon } X_{\theta(R_n)} \right) \\
& = \widehat{\mathbf{E}}^\uparrow  \left( \widehat{\mathbf{P}}^\uparrow  \left( J_{R_{n}} \leq \frac{\epsilon}{1- \epsilon } X_{\theta(R_n)}\Big\vert X_{\theta(R_n)}\right) \right) \\
&=  1 -\widehat{\mathbf{E}}^\uparrow \left( \frac{W(\ell(\epsilon) X_{\theta(R_n)})}{W( X_{\theta(R_n)})} \right), \\
\end{split}
\]
 where $\ell(\epsilon) = (1- 2 \epsilon)/(1- \epsilon)$.
Since the hypothesis $(\mathbf{H})$ is satisfied, an application
of Fatou-Lebesgue Theorem shows that
\begin{eqnarray*}
 \limsup_{n \rightarrow + \infty}\widehat{\mathbf{E}}^\uparrow \left( \frac{W(\ell(\epsilon) X_{\theta(R_n)})}{W( X_{\theta(R_n)})} \right) \leq  \widehat{\mathbf{E}}^\uparrow \left(\limsup_{n\to\infty} \frac{W(\ell(\epsilon) X_{\theta(R_n)})}{W( X_{\theta(R_n)})} \right) < 1,
\end{eqnarray*}
 which implies that
\begin{eqnarray*}
  \lim_{n \rightarrow \infty}\widehat{\mathbf{P}}^\uparrow \left( \frac{X_{\theta(R_n)}-J_{R_{n}}}{f(R_{n})} \geq (1 - 2\epsilon) \right) > 0.
\end{eqnarray*}
Next, we note
\begin{eqnarray*}
 \widehat{\mathbf{P}}^\uparrow \left( \frac{X_{\theta(R_p)} -J_{R_{p}}}{f(R_{p})} \geq (1 - 2 \epsilon),  \mbox{ for some } p \geq n \right) \geq \widehat{\mathbf{P}}^\uparrow \left( \frac{X_{\theta(R_n)}-J_{R_{n}}}{f(R_{n})} \geq (1 - 2\epsilon) \right).
\end{eqnarray*}
Since $R_{n}$ converges to 0, $\widehat{\mathbf{P}}^\uparrow$-a.s. as $n$ goes to $\infty$,
it is enough to take limits in both sides. Therefore for all $\epsilon \in (0,1/2)$
\[
\widehat{\mathbf{P}}^\uparrow \left( \frac{X_{\theta(t)}-J_{t}}{f(t)}  \geq  (1 - 2 \epsilon), \hspace*{2mm} \mbox{i.0.,} \hspace*{2mm} \mbox{as} \hspace*{2mm} t \rightarrow 0 \right) \geq \lim_{n \rightarrow + \infty}  \widehat{\mathbf{P}}^\uparrow \left( \frac{X_{\theta(R_n)}-J_{R_{n}}}{f(R_{n})}  \geq (1 - 2\epsilon) \right)=1.\\
\]
Then,
\[ \limsup_{t \rightarrow 0} \frac{X_{\theta(t)}- J_{t}}{f(t)} \geq 1-2\epsilon, \hspace*{4mm} \widehat{\mathbf{P}}^\uparrow -\mbox{a.s.,} \]
 and since $\epsilon$ can be chosen arbitrarily small,
we get the result.
\hfill $\Box$

\section{Concluding remarks on quasi-stationarity}

We conclude this paper a  brief remark about a kind of
conditioning of CB-process which result in a so-called
quasi-stationary distribution. Specifically we are interested in
establishing the existence of a normalization constant $\{ c_{t},
t\geq 0 \}$ such that the weak limit
\[
 \lim_{t \to \infty} \mathbb{P}_{x}( Y_{t}/c_{t} \in dz |
T_{0} > t ),
 \]
exist for $x >0$ and $z \geq 0$. 

Results of this kind have been established for CB-processes for
which the underlying spectrally positive Le\'vy process has a second
moment in  \cite{lamb}; see also \cite{li}, and for the $\alpha$-stable CB-process
 with $\alpha \in (1,2]$ in \cite{kyp}. In the more
general setting, \cite{pakes} formulates conditions for the
existence of such a limit and characterizes the resulting
quasi-stationary distribution. The result below shows, when the
branching mechanism is regularly varying at $\infty$, an explicit
formulation of the normalization sequence $\{c_{t}: t \geq 0\}$ and
the limiting distribution is possible.

\begin{lemma}
Suppose that the branching mechanism $\psi$ is regularly varying at
 $\infty$ with index $\alpha \in (1,2]$. Then, for all $x \geq 0$, with $c_t= 1/ \varphi(t)$
\[
 \lim_{t\to\infty}\mathbb{E}_x\Big[e^{-\lambda Y_t/c_t}\Big|T_0>t\Big] = 1 - \frac{1}{[1+ \lambda^{-(\alpha-1)}]^{1/(\alpha-1)}}.
\]
\end{lemma}

\noindent \emph{Proof:} The proof pursues a similar line of
reasoning to the the aforementioned references \cite{kyp, lamb, li,
pakes}. From (\ref{CBut}) it is straightforward to deduce that
\[
 \lim_{t\to\infty}\mathbb{E}_x\Big[1- e^{-\lambda Y_t/c_t}\Big|T_0>t\Big] = \lim_{t\to\infty}\frac{u_t(\lambda/c_t)}{u_t(\infty)} = \lim_{t\to\infty}\frac{\varphi(t + \phi(\lambda/c_{t}))}{\varphi(t)},
\]
if the limit on the right hand side exists. However, since $\psi$ is
regularly varying at $\infty$ with index $\alpha$, we have that
$1/\psi$ is regularly varying at $\infty$ with index $-\alpha$. On
the other hand an application of Karamata's Theorem (see for instance
Bingham et al \cite{bi2}) gives
\[ 
\phi(t) \sim - \frac{1}{1-\alpha}\frac{t}{\psi(t)} \hspace*{4mm}
\mbox{as} \hspace*{4mm} t \rightarrow \infty,
 \]
 and therefore,
\[ 
\lim_{t \rightarrow \infty} \frac{\phi(\lambda t)}{\phi(t)} =
\lambda^{1-\alpha} \hspace*{4mm} \mbox{for} \hspace*{3mm} \lambda > 0.
 \]
In other words $\phi$ is regularly varying at $\infty$ with index
$1-\alpha$ and then its inverse $\varphi$ is regularly varying at
$\infty$ with index $\frac{1}{1 - \alpha}$, thus for all $\lambda, \epsilon
>0$,   and  $t$ sufficiently
large, we have
\[
 \left(\frac{\lambda}{1-\epsilon} \right)^{-(\alpha-1)} t \leq \phi( \lambda \varphi(t)) \leq \left(\frac{\lambda}{1+\epsilon} \right)^{-(\alpha-1)} t. 
 \]
Therefore since $\varphi$ is decreasing, we deduce
\[
\frac{\varphi \left( \left(1 + \left( \lambda/(1+\epsilon) \right)^{-(\alpha-1)} \right) t   \right)}{ \varphi(t)} \leq \frac{\varphi(t + \phi(\lambda/c_{t}))}{\varphi(t)} \leq \frac{\varphi \left( \left(1 + \left( \lambda/(1-\epsilon) \right)^{-(\alpha-1)} \right) t   \right)}{ \varphi(t)},
\]
for $t$ sufficiently large. On the other hand,  since $\varphi$ is regular
varying we deduce
\[
\frac{1}{ \left(1 + \left( \frac{\lambda}{1+\epsilon} \right)^{-(\alpha-1)} \right)^{\frac{1}{\alpha-1}}  } \leq \frac{\varphi(t + \phi(\lambda/c_{t}))}{\varphi(t)} \leq \frac{1 }{ \left(1 + \left( \frac{\lambda}{1-\epsilon} \right)^{-(\alpha-1)} \right)^{\frac{1}{\alpha-1}}  },
\]
\noindent for  $t$
sufficiently large.  Hence the result follows taking the limit in the above inequality as $\epsilon$ goes to $0$.
\hfill $\Box$ \\

\end{document}